\documentclass[10pt,a4paper]{article}
\usepackage{amsfonts}
\usepackage[centertags]{amsmath}
\usepackage{amsfonts}
\usepackage{amssymb}
\usepackage{amsthm}
\usepackage{newlfont}

\newcommand{\comment}[1]{}

\newcommand{\zd}{\mathbb{Z}^{d}}
\newcommand{\zii}{\mathbb{Z}^{2}}
\newcommand{\ld}{\mathbb{L}^{d}}
\newcommand{\ed}{\mathbb{E}^{d}}

\newcommand{\la}{\Lambda}

\newcommand{\E}{\mathbb{E}}

\newcommand{\F}{\mathcal{F}}

\newcommand{\mlo}{\mu_{\Lambda_{o}}}
\newcommand{\cyi}{C_{Y}^{+}}

\newcommand{\bi}{B^{+}}
\newcommand{\bii}{B^{-}}

\newtheorem{dfn}{Definition}[section]
\newtheorem{teo}[dfn]{Theorem}
\newtheorem{prop}[dfn]{Proposition}
\newtheorem{lem}[dfn]{Lemma}
\newtheorem{rem}[dfn]{Remark}
\newtheorem{problem}[dfn]{Problem}

\begin{document}

\title{{\LARGE\sf  Slab Percolation and Phase Transitions for the Ising
 Model
}}

\medskip

\author{\textbf{Emilio De Santis, Rossella Micieli} \\{\small
\textsl{Dipartimento di Matematica ``Guido Castelnuovo'',} }
\\{
\small \textsl{ Universit\`{a} di Roma} \textit{La Sapienza\/}, }
\\{\small \textsl{Piazzale Aldo Moro, 2}}
{\small \textsl{- 00185 Roma, Italia}}
\\{\small \textsl{ e-mail: desantis@mat.uniroma1.it  }}  }

\date{June 8, 2004} 
\maketitle

\begin{abstract}
We prove, using the random-cluster model, a strict inequality
between site percolation and magnetization in the region of phase
transition for the $d$-dimensional Ising model, thus improving a
result of \cite{cnpr}. We extend this result also at the case of
two plane lattices $\zii$ (slabs) and give a characterization of
phase transition in this case.  The general case of $N$ slabs,
with $N$ an arbitrary positive integer, is partially solved and it
is used to show that this characterization holds in the case of
three slabs with periodic boundary conditions. However in this
case we do not obtain useful inequalities between magnetization
and percolation probability.
\end{abstract}

\maketitle
\section{Introduction}

At the end of the 70's the seminal paper \cite{cnpr} showed the
connection between phase transition for the ferromagnetic Ising
model and site percolation. This point of view has given a
geometrical interpretation of phase transition, initiating a new
line of research. Following this approach Higuchi developed
techniques to study percolation for the two dimensional Ising
model, with non zero external field \cite{Hig0, Hig1,Hig2, Hig3}.
For example, he showed that for every $\beta<\beta_{c}$ there
exists a positive critical point $h_c(\beta)$ such that an
infinite cluster of (+)-sites does not exist for all $h<h_c$.
Russo \cite{Ru}, in 1979, proved that if an equilibrium Ising
measure is invariant under translations along one direction of the
two dimensional lattice then it is invariant under all
translations. Hence, Aizenman and Higuchi (see \cite{Ai, Hignn}
and also the more recent paper \cite{GH}) showed that the only
extremal Gibbs measures are $\mu_+$ and $\mu_-$. In three
dimensions the situation is different: Dobrushin showed that at
low temperatures there exist non-translation invariant Gibbs
states.  There are substantial differences between two and three
dimensions also for percolation for the Ising measure; in fact in
\cite{CR} it is showed that in three dimensions there is
coexistence of infinite plus and minus clusters (at least for
small values of the parameter $\beta$), while in two dimensions it
is proved that the infinite clusters of opposite sign can not
coexist \cite{cnpr, non.coesistenza}.

\medskip

The paper is organized as follows: in Section \ref{basic} we set
the notation and introduce some basic tools.  In Section \ref{zd}
we use the random cluster model to give, below the critical
temperature, a strict inequality between magnetization and site
percolation probability for the $d$-dimensional cubic lattice, in
this way we improve a result of \cite{cnpr}.  Then, in Section
\ref{Ns} we partially generalize the result to some \emph{slab
graphs}. A slab graph is a graph $G_N=(V_N,E_N)$, where
$V_N=\zii\times\{0,\ldots,N-1\}$ and $E_N$ is the set of all pairs
of vertices in $V_{N}$ having Euclidean distance equal to one.
Such graphs, also called \emph{bunkbed graphs}, have attracted the
attention of other researches (see \cite{Boo, HagB}) in the study
of random walks, random-cluster model and some correlation
inequalities for the Ising model. In Section~\ref{ultima}, for $N
= 2 $, we prove a characterization of phase transitions similar to
that for the ${\mathbb{ Z}^2}$ lattice, thereby obtaining an
inequality between magnetization and percolation probability of
{\it columns} formed only of +1 spins. For $N =3$ and periodic
boundary conditions, we are able to characterize the phase
transition through percolation of columns with majority of plus.
However we cannot obtain any meaningful inequality in this case.

Recently, Bodineau has proved a fine and natural result on slab
percolation for the Ising model (see \cite{Bo}); let
$\beta_{c}^{(N)} $ be the critical point for the graph $G_N $ and
$\beta_{c} (3)$ the critical point for the three dimensional cubic
lattice. It is easy to show, using the FKG inequality for the
random cluster measure, that $ \beta_{c}^{(N)} >
\beta_{c}^{(N+1)}>\cdots \geq \beta_{c} (3) $ for every $N $.
However Bodineau also show that $\lim_{N \to \infty}
\beta_c^{(N)}=\beta_{c} (3) $; thus the $N$ slabs are a good
approximation for the three dimensional Ising model, at least for
the purpose of estimating the critical point.

\comment{

The study of phase transitions in the Ising model is simplified by
the utilization of the random-cluster model.

In section 2 we give some basic definitions and the notation used
in the next sections.

In section 3 we prove a result given by Coniglio, Nappi, Peruggi,
Russo \cite{cnpr} in an other way, by using the random-cluster
model. The proof is really straightforward. Moreover we show that
for $T<T_{c}(d)$ the percolation probability is strictly greater
than magnetization. If $d=2$ we have a complete characterization
of phase transition in the Ising model.

The following sections are dedicated to slabs results. The
percolation probability, referred to slabs, is defined as the
probability that the ``column'' at the origin (where a ``column''
is a pair of vertices with the same plane coordinates) belongs to
an infinite cluster of ``columns'' with a majority of vertices
with spin +1 (or $-1$).

In section 4 we extended the characterization of phase transition
at the case of two slabs. In section 5 we analyze the case of $N$
slabs and give some results. Finally, in section 6, by adding an
hypothesis, we prove the characterization of phase transition in
the case of three slabs using the results found in section 5.

}
\section{Basic definitions and notation}   \label{basic}

In this section we set our notation for percolation, ferromagnetic
Ising model, and random-cluster model.

Let $d\geq 2$ and let $\zd$ be the set of all $d$-vectors
$x=(x_{1},x_{2},\ldots,x_{d})$ with integral coordinates. The
distance $||x-y||$ from $x$ to $y$ is defined by $
||x-y||=\sum_{i=1}^{d}|x_{i}-y_{i}|$.  If $||x-y||=1$ we say that
$x$ and $y$ are \emph{adjacent}. We turn $\zd$ into a graph,
called the \emph{$d$-dimensional cubic lattice}, by adding edges
$e=\langle x,y\rangle$ between all pairs $x,y$ of adjacent points
of $\zd$; we denote this lattice by $\ld=(\zd,\ed)$, where $\ed$
is the edge set.
The edge $e= \langle x,y\rangle $ is said to be \emph{incident} to
the vertices $x$ and $y$; in this case we also say that $x$ and $y
$ are endvertices of the edge $e \in \ed $.


\medskip \medskip

A \emph{path} of $\ld$ is an alternating sequence
$x_{0},e_{0},\ldots,e_{n-1},x_{n}$ of distinct vertices  and edges
with $e_{i}=\langle x_{i},x_{i+1}\rangle$ for all
$i=0,\ldots,n-1$; such a path has \emph{length n} and is said to
\emph{connect} $x_{0}$ to $x_{n}$. A subset $Y \subset\zd$ is
\emph{connected} if for all pairs $x,y$ of vertices in $Y$, there
exists a path connecting the vertices $x,y$ having all its
vertices belonging to $Y$. The \emph{boundary} of $Y\subset\zd$ is
the set $\partial Y$ of all vertices in $\zd\setminus Y$ that are
adjacent to at least one vertex in $Y$.

\medskip \medskip

The {\it edge configuration space } is $\Omega=\{0,1\}^{\E^{d}}$,
so its elements are vectors $\omega=(\omega(e):e\in\ed)$.  We say
that the edge $e$ is \emph{open} if $\omega(e)=1$, and
\emph{closed} if $\omega(e)=0$. For $\omega\in\Omega$, we consider
the random subgraph of $\ld$ containing the vertex set $\zd$ and
the open edges ($\omega^{-1}(1)$) only; an \emph{open cluster} of
$\omega$ is a maximal connected component of this graph.

Let $\Sigma=\{-1,+1\}^{\zd}$ be the { \it spin configuration
space}, elements of which are $\sigma=(\sigma_{x}:x\in\zd)$. We
say that the vertex $x$ has spin $+1$ ($-1$) if $\sigma_{x}=+1$
($-1$) . For $\sigma\in\Sigma$, we consider the random subgraph of
$\ld$ containing the edge set $\ed$ and the vertices
$\sigma^{-1}(+1)$ only; a \emph{$(+)$-cluster} of $\sigma$ is a
maximal connected component of this graph. A \emph{$(-)$-cluster}
is defined in a similar way.  We use the notation
$(\infty,\pm)$-cluster to indicate an infinite $(\pm)$-cluster.

The spaces $\Omega$ and $\Sigma$ are endowed with the discrete
topology. We denote by $\F$ the $\sigma$-field generated by the
finite-dimensional cylinders. For $A\in\F$, we indicate with
$\bar{A}$ the complement of $A$.

\medskip

Let $\la$ be a finite box of $\zd$, {\it i.e.}
$\la=\prod_{i=1}^{d}[x_{i},y_{i}]$ for some $x,y\in\zd$, where
$[x_{i},y_{i}]$ is the set
$\{x_{i},x_{i}+1,x_{i}+2,\ldots,y_{i}\}$. We write $\E_{\la}$ for
the set of edges $e=\langle x,y\rangle$ in $\ed$ such that
$x,y\in\la$ and we define
\begin{equation}
 \Omega_{\la}^{1}=\{\omega\in\Omega:\omega(e)=1\:\:\textrm{for
all}\:\:e\in\ed\setminus\E_{\la}\},
\end{equation}
\begin{equation}
\Sigma_{\la}^{+}=\{\sigma\in\Sigma:\sigma_{x}=+1\:\:\textrm{for
all}\:\:x\in\zd\setminus\la\}.
\end{equation}
$\Sigma_{\la}^{-} $ is defined analogously.

 For $0\leq p\leq 1$, let $\phi_{\la,p}^{1}$ be the random-cluster
measure on $\Omega_{\la}^{1}$ with wired boundary conditions
\cite{grimmett}, and let $\mu_{\la,\beta,J}^{\pm}$ (or simply
$\mu_{\la}^{\pm}$) be the Ising Gibbs measure on
$\Sigma_{\la}^{\pm}$ with $(\pm)$-boundary conditions, zero
external field ($h=0$) and interactions $\{J_e\}_{e \in {\mathbb
E}_{\la\cup\partial\la}}$\cite{ligget}. In this paper we agree
that on each edge $e$ of the graph under consideration there is a
constant interaction $J_e \equiv J = 1$. In some cases one could
take different values of the interactions on different edges; this
will be partially discussed in the last section.

For $p=1-\exp\{-2\beta \} $, we define the {\it coupling between
Ising and random-cluster measures} $\nu_{\la,p}^{+}$ (or simply
$\nu_{\la}^{+}$) on $\Sigma_{\la}^{+}\times\Omega_{\la}^{1}$ with
$(+)$-boundary conditions, as
\begin{equation}\label{coupl}
\nu_{\la}^{+}(\sigma,\,\omega)=K\,
\prod_{e\in\E_{\la\cup\partial\la}}\{(1-p)\delta_{\omega(e),0}+p\,
\delta_{\omega(e),1}\delta_{e}(\sigma)\},
\end{equation}
where
$K$ is the normalizing constant, and
$\delta_{e}(\sigma)=\delta_{\sigma_{x},\sigma_{y}}$ for $e=\langle
x,y\rangle$ ($\delta_{i,j}$ is the Kronecker delta).
Similarly, we define the coupling measure $\nu_{\la}^{-}$ with
$(-)$-boundary conditions.  The marginal measure of
$\nu_{\la}^{+}$ on $\Sigma_{\la}^{+}$ is the Ising measure
$\mu_{\la}^{+}$, and the marginal measure on $\Omega_{\la}^{1}$ is
the random-cluster measure $\phi_{\la,p}^{1}$. We say that an open
cluster of $\omega\in\Omega_{\la}^{1}$ (or a ($+$)-cluster of
$\sigma\in\Sigma_{\la}^{+}$) \emph{touches} $\partial \la$ if at
least one vertex of this cluster belongs to $\partial\la$. Given
$\omega$, the conditional measure of $\nu_{\la}^{+}$ on
$\Sigma_{\la}^{+}$ is obtained by putting $\sigma_{x}=+1$ for
every $x\in\partial\la$, then by setting $\sigma_{x}=+1$ for every
$x$ in an open cluster touching $\partial\la$, and finally by
choosing spins $+1$ or $-1$ with probability $\frac{1}{2}$ on all
open clusters not touching $\partial\la$. Given $\sigma$, the
conditional measure of $\nu_{\la}^{+}$ on $\Omega_{\la}^{1}$ is
obtained by setting $\omega(e)=1$ for every
$e\in\E_{\partial\la}$, then for $e\in\E_{\la}$, $\omega(e)=0$ if
$\delta_{e}(\sigma)=0$ and $\omega(e)=1$ with probability $p$
(independently of other edges) if $\delta_{e}(\sigma)=1$.

Using standard arguments on the stochastic order (FKG order) one
can define these measures directly on the infinite spaces by
taking the weak limit of measures; thus are well defined: $ (
\Sigma \times \Omega , {\mathcal F}, \nu_+ ) $, $ ( \Omega ,
{\mathcal F}, \phi^1_p ) $ and $ ( \Sigma , {\mathcal F}, \mu_+ )
$, where $\phi^1_p$, $\mu_+$ and $\nu_+$ are respectively the
random-cluster measure on $\Omega$ with wired boundary conditions,
the Ising measure on $\Sigma$ with $(+)$-boundary conditions, and
the coupling measure on  $\Sigma \times \Omega$ with
$(+)$-boundary conditions.  Moreover for any given $\omega \in
\Omega $, the conditional probability measure $\nu_+ (\cdot |
\omega ) $ is obtained by setting $\sigma_x = +1 $ for all the
vertices $x$ belonging to an infinite open cluster of $\omega$,
and by putting spins $+ 1 $ or $-1$  with probability
$\frac{1}{2}$ on each other cluster.

For more details and in the more general setting of every boundary
condition and no-ferromagnetic interactions see \cite{newman} in
which all constructions and relations between the three measures
$\nu_+$, $\mu_+$ and $\phi^1_p$ are clearly explained.

We define $\tau_{\la}=\mu_{\la}^{+}(\sigma_{0}=+1)-\frac{1}{2}$,
where $0$ denotes the origin of $\zd$. If
$\tau=\lim_{\la\uparrow\zd}\tau_{\la}$ (the limit  exists by
monotonicity) then
\begin{equation}\label{tau}
\tau =  \mu_{+}(\sigma_{0}=+1)-\frac{1}{2}.
\end{equation}

Given $X$, $Y$ subsets of $\zd$ we denote with $\{X\leftrightarrow
Y\}$ the set of configurations $\omega\in\Omega$ such that there
exists a vertex $x\in X$ connected to a vertex $y\in Y$ by a path
of open edges. We write $\{x\leftrightarrow\infty\}$ for the set
of configurations $\omega\in\Omega$ such that $x$ belongs to an
infinite open cluster. The following Proposition is a result due
to Kasteleyn and Fortuin (see \cite{grimmett}).
\begin{prop}\label{fk} If $p=1-\exp\{-2\beta \}$, then
$ \tau_{\la}=\frac{1}{2}\,\phi_{\la,p
}^{1}(0\leftrightarrow\partial\la) $.  Equality holds also in the
limit $\la\uparrow\zd$: $
\tau=\displaystyle{\frac{1}{2}}\,\phi_{p}^{1}(0\leftrightarrow\infty).
$
\end{prop}


\medskip
We give some other definitions. We put $
C_{\infty}^{\pm}=\{\sigma\in\Sigma:0\in
\:\textrm{$(\infty,\pm)$-cluster of}\: \sigma\} $.  The
\emph{percolation probability} is denoted by
$R(\pm;\mu_{\pm})=\mu_{\pm}(C_{\infty}^{\pm})$, and the
\emph{magnetization in the origin} is
\begin{equation}\label{magnbis}
M(\mu_{\pm})= {\mathbb E_{ \mu_{\pm} } } (\sigma_0) =
\mu_{\pm}(\sigma_{0}=+1)-\mu_{\pm}(\sigma_{0}=-1),
\end{equation}
where $\sigma_0 $ is the spin on the origin. By (\ref{tau}),
(\ref{magnbis}) and Proposition~\ref{fk} follows
\begin{equation}\label{mt}
M(\mu_{+})=2\tau=\phi_{p}^{1}(0\leftrightarrow\infty).
\end{equation}

\section{$\zd$ percolation and magnetization}\label{zd}
In this section we use the random-cluster model to prove, in a
different way, the inequality relating percolation probability and
magnetization in the $d$-dimensional Ising model, given in
\cite{cnpr}.  Moreover, we prove that for $T<T_{c}(d)$ this
inequality is strict. For $d=2$ we have a complete
characterization of phase transition in the Ising model through
percolation.
\begin{teo}\label{teoc}
 For a ferromagnetic Ising model on $(\Sigma,\F,\mu_{\pm})$ with zero
  external field, the following inequality holds:
\begin{equation}
|M(\mu_{\pm})|\leq R(\pm;\,\mu_{\pm}).
\end{equation}
\end{teo}
\begin{proof}
We consider the following event:
\begin{equation}\label{incl}
  \Sigma \times \{\omega
  \in \Omega:0\leftrightarrow \infty \}
  =C^{+}_{\infty}\times \{\omega
  \in \Omega:0\leftrightarrow \infty \} \cup
 \overline{C_{\infty}^{+}}  \times \{\omega
  \in \Omega:0\leftrightarrow \infty \}
\end{equation}
and we observe that
\begin{equation}\label{zero}
\nu_{+} ( \overline{C_{\infty}^{+}}  \times
\{0\leftrightarrow\infty \})= 0.
\end{equation}
Therefore
\begin{equation}\label{coupling}
\phi^1_p (0\leftrightarrow\infty   ) =    \nu_{+}(\Sigma \times
\{0\leftrightarrow\infty \} ) \leq\nu_{+}(C_{\infty}^{+} \times
\Omega ) = \mu_+ ( C_{\infty}^{+}  )  .
\end{equation}
\comment{
 But on the
other hand,
\begin{equation}
\nu_{\la}^{+}(\Sigma^+_{\Lambda} \times
\{0\leftrightarrow\partial\la
\})=\phi_{\la,p}^{1}(0\leftrightarrow\partial\la)
\end{equation}
and
\begin{equation}
\nu_{\la}^{+}(C_{\partial\la}^{+} \times \Omega_{\Lambda}^1
)=\mu_{\la}^{+}(C_{\partial\la}^{+}),
\end{equation}
so, in the limit $\la\to\zd$, we obtain by (\ref{coupling})
\begin{equation}\label{dis}
\phi_{p}^{1}(0\leftrightarrow\infty)\leq\mu_{+}(C_{\infty}^{+})=R(+;\mu_{+}).
\end{equation}

} Thus, by (\ref{mt}) and (\ref{coupling}) follows
\begin{equation}
M(\mu_{+})=\phi_{p}^{1}(0\leftrightarrow\infty)\leq \mu_+(
C_{\infty}^{+}  )   =   R(+;\mu_{+}).
\end{equation}
Similarly, taking $(-)$-boundary conditions we have $
|M(\mu_{-})|\leq R(-;\,\mu_{-})$.
\end{proof}

 Now, we prove that below the critical temperature the percolation
 probability is strictly larger than the magnetization, and we give a
 characterization of phase transition for $d=2$.
\comment{
 We remark that a
 probability measure $\mu$ on ($\Sigma,\F)$ is \emph{ergodic} if for
 all $W\in\F$, which is invariant under translation, then $\mu(W)=0$
 or $1$.
}
\begin{teo}\label{one}
For a ferromagnetic Ising model on $(\Sigma,\F,\mu_{\pm})$, at
zero external field, the following relations hold:\begin{itemize}
\item[(i)]
$R(\pm;\,\mu_{\pm})\geq|M(\mu_{\pm})|+\displaystyle{\frac{1}{2}}
|M(\mu_{\pm})|\left(\frac{p}{2-p}\right)^{2d(3^{d}-1)}(1-p)^{2d}$,\\
where $p=1-\exp\{-2\beta \}$; \item [(ii)] if $d=2$, then
\begin{equation*}
R(\pm;\,\mu_{\pm})>0\Leftrightarrow|M(\mu_{\pm})|>0.
\end{equation*}
\end{itemize}
\end{teo}
\begin{proof} We prove claim (i) first.  Let
\begin{equation*}
\la'=\{x\in\zd\setminus\{0\}:-1\leq x_{i}\leq 1\:\:\textrm{for all
$i=1,\ldots,d$}\}.
\end{equation*}
 We consider the following cylinders on $\la'$
\begin{eqnarray*}
& & \!\!\!\!\!\!\!\!A=\{\omega\in\Omega:\omega(e)=1\:\:\textrm{for
$e\in\E_{\la'}$},\,\omega(e)=0\:\:\textrm{for $e=\langle
0,y\rangle,y\in\la'$}\},\\
& & \!\!\!\!\!\!\!\!B=\{\omega\in\Omega:\omega(e)=1\:\:\textrm{for
$e\in\E_{\la'}$},\,\omega(e)=1\:\:\textrm{for $e=\langle
0,y\rangle,y\in\la'$}\}.
\end{eqnarray*}
Note that the event $A$ (resp. $B$)  forces  all edges in $\la'$
to be open (resp. open), and all edges incidents at the origin to
be closed (resp. open).

 Let $x_{0}$ be a vertex adjacent to the origin
\setlength\arraycolsep{2pt}
\begin{eqnarray*}
\{\sigma_{0}=+1\}\times(\{x_{0}\leftrightarrow\infty \}\cap A)&=&
[(C_{\infty}^{+}\cap\{\sigma_{0}=+1\})\times
(\{x_{0}\leftrightarrow\infty\}\cap A)]\\
& & \cup[(\overline{C_{\infty}^{+}}\cap\{\sigma_{0}=+1\})\times
(\{x_{0}\leftrightarrow\infty\}\cap A)],
\end{eqnarray*}
and it is clear that
\begin{equation*}
\nu_{+}( (\overline{C_{\infty}^{+}}\cap\{\sigma_{0}=+1\})\times
(\{x_{0}\leftrightarrow\infty\}\cap A) )=0 .
\end{equation*}
 Thus, by using (\ref{incl}) and
noting that events
$\{\sigma_{0}=+1\}\times(\{x_{0}\leftrightarrow\infty\}\cap A)$
and $\Sigma\times\{0\leftrightarrow\infty \}$ are disjoint, we
obtain
\begin{equation*}
\nu_{+}(C_{\infty}^{+}\times\Omega)\geq\nu_{+} (\Sigma\times
\{0\leftrightarrow\infty \})+
\frac{1}{2}\nu_{+}(\Sigma\times(\{x_{0}\leftrightarrow\infty\}\cap
A)),
\end{equation*}
hence,\setlength\arraycolsep{2pt}
\begin{eqnarray}\label{tesis}
R(+;\mu_{+})=\mu_{+}(C_{\infty}^{+})&\geq&\phi_{p}^{1}
(0\leftrightarrow\infty)+\frac{1}{2}\phi_{p}^{1}
(\{x_{0}\leftrightarrow\infty\}\cap A)=\\ &=&
M(\mu_{+})+\frac{1}{2}\phi_{p}^{1}(\{x_{0}\leftrightarrow\infty\}\cap
A).\nonumber
\end{eqnarray}
For the structure of the random cluster measure we obtain
\begin{equation}\label{1}
\phi_{p}^{1}(x_{0}\leftrightarrow\infty\,|\,
A)=\phi_{p}^{1}(x_{0}\leftrightarrow\infty\,|\, B).
\end{equation}
Events $\{x_{0}\leftrightarrow\infty\}$ and $B$ are increasing,
thus by FKG inequality \cite{FKG} we obtain
\begin{equation}\label{2}
\phi_{p}^{1}(x_{0}\leftrightarrow\infty\,|\,
B)\geq\phi_{p}^{1}(x_{0}\leftrightarrow\infty).
\end{equation}
By (\ref{1}) and (\ref{2}) follows
\begin{equation}\label{3}
\phi_{p}^{1}(\{x_{0}\leftrightarrow\infty\}\cap
A)\geq\phi_{p}^{1}(x_{0}\leftrightarrow\infty)\,
\phi_{p}^{1}(A)=M(\mu_{+})\phi_{p}^{1}(A),
\end{equation}
where the last equality follows by the translation invariance of
$\phi_{p}^{1}$.

 We prove now a lower bound for $\phi_{p}^{1}(A)$. Let $k(d,\la')$ the
number of edges of the graph $(\la',\E_{\la'})$. For
$\omega_{\setminus e}\in\Omega_{\setminus
e}=\{0,1\}^{\ed\setminus\{e\}}$, we have \cite{for} (see also
\cite{grimmett, newman})
\begin{equation}\label{4}
\phi_{p}^{1}(\omega(e)=1\,|\,\omega_{\setminus
e})\in\left\{p,\frac{p}{2-p}\right\}.
\end{equation}
Moreover
\begin{equation}\label{5}
 \phi_{p}^{1}(A)=\phi_{p}^{1}(\{\omega(e)=1,e\in\E_{\la'}\}
\cap\{\omega(e)=0,e=\langle 0,y\rangle,y\in\la'\}),
\end{equation}
thus, by (\ref{4}) and (\ref{5}) follows
\begin{equation}\label{6}
\phi_{p}^{1}(A)\geq
 \left(\frac{p}{2-p}\right)^{k(d,\la')}(1-p)^{2d} .
\end{equation}
By (\ref{3}) and (\ref{6}),  we obtain
\begin{equation}\label{tesii}
\phi_{p}^{1}(\{x_{0}\leftrightarrow\infty\}\cap A)\geq
M(\mu_{+})\left(\frac{p}{2-p}\right)^{k(d,\la')}(1-p)^{2d}.
\end{equation}

We also give an upper bound for $k(d,\la')$. The number of
vertices in $\la'$ is $3^{d}-1$ and there are at most $2d$ edges
incident to each vertex in $\la'$, so $k(d,\la')\leq 2d(3^{d}-1)$.
Thus, by (\ref{tesis}) and (\ref{tesii}) follows (i) for
$(+)$-boundary conditions. In a similar way (i) can be proved for
$(-)$-boundary conditions.

We can now prove claim (ii).
If $|M(\mu_{\pm})|>0$ then $R(\pm;\,\mu_{\pm})>0$ by (i).
Conversely, we assume $M(\mu_{+})=M(\mu_{-})=0$ and prove that
$R(+;\,\mu_{+})=R(-;\,\mu_{-})=0$. If $M(\mu_{+})=0$, then
$\mu_{+}=\mu_{-}=\mu$ because there is not phase transition(see
\cite{ligget}). Suppose that $R(+;\mu_{+})=R(+;\mu)>0$, hence also
$R(-;\mu)>0$. Then
\begin{equation*}
 \mu(\exists
\hbox{ 
 $(\infty,\pm)$-cluster}) \geq R(\pm;\mu)>0 .
\end{equation*}
As the events $\{ \exists \hbox{ } (\infty,+)$-cluster$\}$ and $\{
\exists \hbox{ } (\infty,-)$-cluster$\}$ are invariants under
translation and $\mu$ is ergodic \cite{ergodicita}, we have
\begin{equation}\label{pi}
\mu(\{ \exists \hbox{ $(\infty, +)$-cluster} \} \cap \{ \exists
\hbox{ $(\infty, -)$-cluster} \} )=1 .
\end{equation}
 On $\zii$, under suitable conditions (see \cite{non.coesistenza}), an
infinite (+)-cluster cannot coexists with an infinite
($-$)-cluster. The conditions are: translation invariance,
ergodicity, FKG inequality and invariance to reflections with
respect to ${\hat x } $, ${\hat y }$ axes.

These conditions are satisfied by the Ising measure $\mu $ (see
\cite{ergodicita}).  This fact contradicts (\ref{pi}), then
$R(+;\mu)=R(-;\mu)=0$. Note that this claim is proved also in
\cite{cnpr2}. We have reported this alternative proof which
immediately follows by the result in \cite{non.coesistenza}.
\end{proof}

Theorem \ref{one} (i) says that if the temperature is lower than
the critical temperature, or equivalently if the magnetization is
positive, then the percolation probability is strictly greater
than the magnetization. Moreover, for $d=2$, Theorem \ref{one}
(ii) gives a characterization of phase transition through
percolation.  We end this section with a
\begin{rem} \label{limitep}
The Onsager solution for the two dimensional ferromagnetic Ising
model show the exact value of magnetization as a function of
$\beta \in [\beta_c, \infty )$ \cite{Ons}. It is
\begin{equation} \label{onsager}
M (\mu_+) = \{1 -[ \sinh (2\beta )  ]^{-4} \}^{\frac{1}{8}} .
\end{equation}
We can re-write (\ref{onsager}) as a function of the parameter $x
= 1-p= \exp\{-2\beta \}$ obtaining
\begin{equation} \label{onsager2}
M (\mu_+) =\left \{1 -\left [ \frac{2 x }{1 - x^2} \right ]^{4}
\right \}^{\frac{1}{8}} .
\end{equation}
Then, using Taylor expansion we obtain $ m =1- 2 x^4 + o (x^4) $,
giving the magnetization for small values of the parameter
temperature (small $x$). We do not have an explicit formula for
the percolation probability but for small $x $ it is easy to
calculate the first terms in  Taylor expansion. We find
 \begin{equation} \label{percp}
R( +; \mu_+ ) = 1 - x^4 + o ( x^4).
\end{equation}
This general relation also holds for regular graphs
$$( 1- R (+ ; \mu_+ )) \sim  1/2 (1-M (\mu_+) )\sim x^n $$
where $n $ is the degree of the origin.
\end{rem}

\section{$N$ slabs percolation and magnetization} \label{Ns}
In this section we propose a conjecture for the characterization
of phase transition through percolation in the case of $N$ slabs
and some partial results.

We introduce some basic definitions for slabs. Let $\zii$ be the
two-dimensional lattice, and consider the set
$\zii\times\{0,1,\ldots,N-1\}$, where $N$ is an arbitrary positive
integer. The set of all the edges with endvertices in
$\zii\times\{0,1,\ldots,N-1\}$ is denoted by
 $\E^{2,N}$.
\begin{dfn}
An \emph{$N$-vertex} $\mathbf{c}_{i,j}$ of
$\zii\times\{0,1,\ldots,N-1\}$ is a vector
\begin{equation*}
\mathbf{c}_{i,j}=((i,j,0),(i,j,1),\ldots,(i,j,N-1))
\end{equation*}
where $i,j\in\mathbb{Z}$. An \emph{$N$-edge} $\mathbf{e}$ of
$\zii\times\{0, 1,\ldots,N-1\}$ is formed by a couple of
N-vertices $< \mathbf{c}_{i,j} , \mathbf{c}_{l,m} >$ where the
vertices $ ( i,j ), \,  (l,m  ) \in {\mathbb Z}^2 $ are adjacent.

\end{dfn}
We put $\Sigma^{(N)}=\{-1,+1\}^{\zii\times\{0,1,\ldots,N-1\}}$,
and for $\sigma\in\Sigma^{(N)}$ we indicate with
$\sigma_{i,j,k}\in\{-1,+1\}$ the spin of the vertex $(i,j,k)$. An
\emph{$N$-path} is an alternating sequence of $N$-vertices and
$N$-edges as in the definition of path but substituting a vertex
with an $N$-vertex and an edge with an $N$-edge. An
\emph{$N$-subset} of $\zii\times\{0,1,\ldots,N-1\}$ is a subset of
$\zii\times\{0,1,\ldots,N-1\}$ formed by $N$-vertices. An
$N$-subset $Y$ of $\zii\times\{0,1,\ldots,N-1\}$ is
\emph{$N$-connected} if for all pairs of $N$-vertices
$\mathbf{c}_{i,j},\mathbf{c}_{r,s}$ in $Y$ there exists an
$N$-path formed by $N$-vertices in $Y$ having
$\mathbf{c}_{i,j},\mathbf{c}_{r,s}$ as terminal $N$-vertices. We
denote with $\Gamma$ the family of all finite $N$-connected
$N$-subset of $\zii\times\{0,1,\ldots,N-1\}$ containing the
$N$-vertex at the origin $\mathbf{c}_{0,0}$. An \emph{$N$-box} is
a finite $N$-subset of $\zii\times\{0,1,\ldots,N-1\}$ as in the
definition of box but substituting a vertex by an $N$-vertex.
Analogously, we define the $N$-boundary of $Y\in\Gamma$, which is
denoted again with $\partial Y$ for simplicity of notation.
\begin{dfn}   \label{alal}
 For
$\sigma\in\Sigma^{(N)}$, a \emph{$(\mathbf{c}^{+})$-cluster}
(\emph{$(\mathbf{c}^{-})$-cluster}) of $\sigma$ is a maximal
 connected component of $N$-vertices $\mathbf{c_{i,j}}$ such as
 \begin{equation*}
\sum_{k=0}^{N-1} \sigma_{i,j,k}>0 \quad\left(\sum_{k=0}^{N-1}
\sigma_{i,j,k}<0\right).
\end{equation*}
We write $(\infty,\mathbf{c}^{\pm})$-cluster for an infinite
$(\mathbf{c}^{\pm})$-cluster.
\end{dfn}

We set
\begin{equation*}
N\hbox{-}C^{\pm}_{\infty}=\{\sigma\in
\Sigma^{(N)}:\mathbf{c}_{0,0}\in \textrm{
$(\infty,\mathbf{c}^{\pm})$- cluster  of $\sigma$}\},
\end{equation*}
where $\mathbf{c}_{0,0}$ is the $N$-vertex at the origin. Notice
that $N\hbox{-}C_{\infty}^{+} $ ($N\hbox{-}C_{\infty}^{-}$) is the
event that the N-origin belongs to an infinite cluster of
$N$-vertices with a majority of spins $+1$ ($-1$) on every
N-vertex. Let $E^{+}$ ($E^{-}$) be the set of configurations in
$\Sigma^{(N)}$such that the $N$-vertex at the origin has a
majority of spins +1 ($-1$) in its vertices. The events $E^{+}$
and $E^{-}$ are disjoint and, for odd values of $N$,  $E^{+}\cup
E^{-}=\Sigma^{(N)}$.

 Let
$\mu_{\pm}$ be the Ising  measure on $\Sigma^{(N)}$ with
$(\pm)$-boundary conditions.  We set up also the {\it vertical
interactions} $J_v \equiv 1 $ between spins on two adjacent vertices
belonging to different slabs.
 The \emph{$N$-percolation probability} is
$R(\mathbf{c}^{\pm};\,\mu_{\pm})=\mu_{\pm}(N-C_{\infty}^{\pm})$.

In  next proposition we show that if the $N$-percolation
probability is positive then  magnetization is positive.

\begin{prop}\label{perc.npiani}
For a ferromagnetic Ising model on $(\Sigma^{(N)},\F,\mu_{\pm})$
at zero external field, the following relation holds:
\begin{equation*}
R(\mathbf{c^{\pm}};\mu_{\pm})>0\Rightarrow |M(\mu_{\pm})|>0.
\end{equation*}
\end{prop}
\begin{proof}
We project the $N$ slabs on a single lattice, $\zii$, by assigning
spins +1 ($-$1) on the vertices corresponding to $N$-vertices with
a majority of spins +1 ($-$1) and choosing spins +1 or $-1$ with
probability $\frac{1}{2}$ on the remaining vertices. This
construction induces a new measure $\pi_{\pm}$ on
$\Sigma=\{-1,+1\}^{\zii}$.  We note that if there exists an
infinite $(\mathbf{c^{+}})$-cluster in $\zii\times\{0,1,\ldots,N-1
\}$, then there exists an infinite $(\pm)$-cluster in the new
lattice.
 If $M(\mu_{\pm})=0$ then $M(\pi_{\pm})=0$. Similarly to Theorem
\ref{one} (ii) by using the result given in \cite{non.coesistenza}
and noting that $\pi_{\pm}$ satisfy all the required hypotheses,
follows $R({\pm};\,\pi_{\pm})=0$.  Thus also
$R(\mathbf{c}^{\pm};\,\mu_{\pm})=0$ by the observation above.

\end{proof}

The opposite implication of Proposition \ref{perc.npiani} will be
partially proved.
\begin{lem}\label{lemma.npiani}
Let $(H,\mathcal{A},\mathbb{P})$ an arbitrary probability space.
If $X$ and $Y$ are random variables  with $X$  symmetric and $Y$
 not negative, then
\begin{equation*}
\mathbb{P}(X+Y>0)\geq\mathbb{P}(X+Y<0).
\end{equation*}
\end{lem}
\begin{proof}
Since $Y\geq 0$, $\{X>0\}\subseteq\{X+Y>0\}$ and
$\{X+Y<0\}\subseteq\{X<0\}$. Thus, because of $X$ is symmetric
\begin{equation*}
\mathbb{P}(X+Y>0)\geq\mathbb{P}(X>0)=\mathbb{P}(X<0)\geq\mathbb{P}(X+Y<0).
\end{equation*}
\end{proof}
The following proposition says that if there is phase transition
then the probability of a majority of $+1$ spins on the $N$-vertex
at the origin (on any fixed $N$-vertex) is larger than the
probability of having a majority of $-1$ spins on such $N$-vertex.
\begin{prop}\label{prop.n+1}
For a ferromagnetic Ising model on $(\Sigma^{(N)},\F,\mu_{\pm})$
at zero external field $|\mu_{\pm}(E^{+}) - \mu_{\pm}(E^{-})|>0$
if and only if $\beta > \beta_{c}(d))$.
\end{prop}
\begin{proof}
Suppose $\beta \leq \beta_{c}(d))$, then $\mu_+=\mu_-=\mu$, so
\begin{equation*}
\mu_\pm(E^+)=\mu(E^+)=\mu(E^-)=\mu_\pm(E^-).
\end{equation*}

Conversely, let us consider $\omega \in
\Omega^{(N)}=\{0,1\}^{\E^{2,N}}$.
 Given $\omega\in\Omega^{(N)}$, the sum of spins on the vertices in
 $\mathbf{c}_{0,0}$ can be expressed as the sum of a symmetric random
 variable (vertices belonging to an finite cluster of open edges) and
 a positive random variable (vertices belonging to an infinite cluster
 of open edges), thus by Lemma \ref{lemma.npiani}
\begin{equation}\label{e+e-}
\nu_{+}(E^{+}\times\Omega^{(N)} \,|\,\omega)\geq
\nu_{+}(E^{-}\times\Omega^{(N)} \,|\,\omega).
\end{equation}
 We have
\begin{equation} \label{cpiu}
 \mu_{+}(E^{+})=\int_{\Omega^{(N)} }
\nu_{+}(E^{+} \times\Omega^{(N)}\,|\,\omega)\, \phi_{p}^{1}( d
\omega),
\end{equation}
\begin{equation}
\mu_{+}(E^{-})=\int_{\Omega^{(N)} } \nu_{+}(E^{-} \times
\Omega^{(N)}\,|\,\omega)\, \phi_{p}^{1}( d \omega) \label{cmeno}
\end{equation}
Let $A$ be the event that all vertices in $\mathbf{c}_{0,0}$
belong to an infinite open cluster. More precisely
\begin{equation*}
A=\{\omega\in\Omega^{(N)}:(0,0,k)\leftrightarrow \infty
\:\:\textrm{for all}\:\:k=0,\ldots,N-1\}.
\end{equation*}
Given $\omega\in A$, the conditional measure is obtained by
setting $\sigma_{0,0,k}=+1$ for every $ k=0,\ldots,N-1$, thus
$\nu_{+}(E^{+} \times\Omega^{(N)}\,|\,\omega)=1$. Hence, by
(\ref{e+e-}) and (\ref{cpiu}), follows \setlength\arraycolsep{2pt}
\begin{equation*}
 \mu_{+}(E^{+})=\int_{A } \nu_{+}(E^{+}
\times\Omega^{(N)}\,|\,\omega)\, \phi_{p}^{1}( d \omega) +
\int_{\Omega^{(N)} \setminus A   } \nu_{+}(E^{+}
\times\Omega^{(N)}\,|\,\omega)\, \phi_{p}^{1}( d \omega) =
\end{equation*}
\begin{equation}\label{prop.fk.i}
 =      \int_{A } \, \phi_{p}^{1}( d \omega) +
\int_{ \Omega^{(N)} \setminus A   } \nu_{+}(E^{+}
\times\Omega^{(N)}\,|\,\omega)\, \phi_{p}^{1}( d \omega)
\geq
\end{equation}
\begin{equation*}
 \geq    \phi_{p}^{1}( A ) + \int_{ \Omega^{(N)} \setminus A  } \nu_{+}(E^{-} \times\Omega^{(N)}\,|\,\omega)\,
\phi_{p}^{1}( d \omega)    .
\end{equation*}
  Moreover if $\omega\in A$ then
$\nu_{+}(E^{-} \times\Omega^{(N)}\,|\,\omega)=0$. So, by
(\ref{cmeno}) follows
\begin{equation}\label{prop.fk.ii}
 \mu_+ (E^-) =    \int_{ \Omega^{(N)} \setminus A }
\nu_{+}(E^{-} \times\Omega^{(N)}\,|\,\omega)\, \phi_{p}^{1}( d
\omega) .
\end{equation}
 By using (\ref{e+e-}), (\ref{prop.fk.i}) and (\ref{prop.fk.ii}), we
obtain
\begin{equation}\label{npianii}
\mu_{+}(E^{+})-\mu_{+}(E^{-})\geq \phi_{p}^{1}(A).
\end{equation}
Consider now the events
\begin{eqnarray*}
& & \!\!\!\!F=\{\omega\in\Omega^{(N)}:(0,0,0)\leftrightarrow\infty \},\\
& &
\!\!\!\!G=\{\omega\in\Omega^{(N)}:\omega(e)=1\:\textrm{for}\:e=\langle
(0,0,k-1),(0,0,k)\rangle,\,k=1,..,N-1\}.
\end{eqnarray*}
We note that
$A\supseteq F\cap G$. Since $F$ and $G$ are increasing events, by
FKG inequality  we obtain
\begin{equation}   \label{disufkg}
\phi_{p}^{1}(A)\geq\phi_{p}^{1}(F\cap
G)\geq\phi_{p}^{1}(F)\,\phi_{p}^{1}(G).
\end{equation}

But, by hypothesis,
$\phi_{p}^{1}(F)=\phi_{p}^{1}((0,0,0)\leftrightarrow\infty)=M(\mu_{+})>0$
and  $\phi_{p}^1(G)>0$ depending on a finite number of edges.  By
 inequality (\ref{disufkg}) we get
$\phi_{p}^{1}(A)>0$, hence
\begin{equation*}
\mu_{+}(E^{+})-\mu_{+}(E^{-})\geq\phi_{p}^{1}(A)>0.
\end{equation*}
The same argument holds for $(-) $-boundary conditions, therefore
\begin{equation*}
\hbox{phase transition} \Leftrightarrow   |M(\mu_{\pm})|>0
\Rightarrow|\mu_{\pm}(E^{+})-\mu_{\pm}(E^{-})|>0.
\end{equation*}
\end{proof}

We are now in the position to present our conjecture for the
characterization
 of phase transition in the Ising model defined on
 the lattices $\zii\times\{0,1,\ldots,N-1\}$. We believe that, for
 these models,
 $N$-percolation probability is positive if and only if there is phase transition.
Proposition \ref{perc.npiani} shows that an implication is true.
To
 prove the other
one we should use Proposition \ref{prop.n+1} and  the next
argument.

 Let $Y\in\Gamma$ be a fixed element of $\Gamma$ and we set
\begin{equation}\label{cy+-}
C^{\pm}_{Y}=\{\sigma\in\Sigma^{(N)}:Y\in\Gamma\:\textrm{is a
$(\mathbf{c}^{\pm})$-cluster of $\sigma$}\}   .
\end{equation}
We have, as in \cite{cnpr}
\begin{eqnarray}\label{e+-}
\mu_{+}(E^{+})-\mu_{+}(E^{-}) & = &
\mu_{+}(E^{+})-\mu_{-}(E^{+})=\\
& = & \sum_{Y\in\Gamma}(\mu_{+}(\cyi)-\mu_{-}(\cyi))+
\mu_{+}(N\hbox{-}C_{\infty}^{+})-\mu_{-}(N\hbox{-}C_{\infty}^{+}).\nonumber
\end{eqnarray}
Thus, a sufficient condition for the claim to hold is that
\begin{equation}\label{impl.npiani}
 \mu_{+}(\cyi)\leq\mu_{-}(\cyi) , \,\,\,\,\,\,\,\,
\,\,\,\,\,\,\,\,\hbox{ for all }  Y \in \Gamma  .
\end{equation}
Indeed, by assumption (\ref{impl.npiani}) and
Proposition~\ref{prop.n+1} for $\beta > \beta_c $ one obtains
\begin{equation}\label{doppio}
R(\mathbf{c}^{+};\,\mu_{+})= \mu_{+}(N\hbox{-}C_{\infty}^{+}) \geq
\mu_{+}(E^{+})-\mu_{+}(E^{-})> 0   .
\end{equation}
Therefore Proposition \ref{perc.npiani} and inequality
(\ref{doppio})  give a characterization of phase transition
through percolation in the case of $N$ slabs. In  next section, we
present the cases of two ($N=2$) and three slabs ($N=3$) with
periodic boundary conditions, showing that (\ref{impl.npiani})
holds.


\section{Two particular cases} \label{ultima}
 In this section, we give a characterization of phase transition
through percolation in the cases of two and three slabs. We manage
to the case of two slabs the result of Theorem~\ref{teoc}, the
proof being similar to \cite{cnpr}. The extension to the case of
three slabs is done in a different flavor.  We start with some new
definitions.

For a fixed $Y\in\Gamma$,  consider an $N$-box $\la_{o}$ such that
$\Lambda_{o}\supset Y\cup\partial Y$. Let
$\Lambda=\Lambda_{o}\cup\partial\Lambda_{o}$ and let
$\mu_{\Lambda}$ be the Ising measure on
$\Sigma^{(N)}_{\Lambda}=\{-1,\,+1\}^{\Lambda}$ with free boundary
conditions. Set
\begin{equation*}
B^{+}=\{\sigma\in\Sigma^{(N)}_{\Lambda}:\sigma_{i,j,0}=+
1,\ldots,\sigma_{i,j,N-1}=+1\:\: \textrm{for all
$\mathbf{c}_{i,j}\in\partial\Lambda_{o}$}\},
\end{equation*}
and similarly for $B^-$. If $\cyi$ is given by (\ref{cy+-}), let
\begin{equation} \label{partt}
\partial\cyi = \{ \sigma \in   \Sigma^{(N)}_{\la} :
\sum_{l= 0}^{N-1} \sigma_{i,j,l} \leq 0  \hbox{ for each }
 {\mathbf c}_{i,j} \in \partial Y   \}
\end{equation}
be the set of configurations in $\Sigma^{(N)}_{\la}$ such that
 each $N$-vertex of $\partial Y$ does not have  a majority of $+1$ spins  on its
 vertices. In general we denote by
\begin{equation*}
\sigma_{V}=\{\tilde{\sigma}\in\Sigma_{\la}^{(N)}:
\tilde{\sigma}_{i,j,k}=\sigma_{i,j,k}\:\textrm{for
all}\:\mathbf{c}_{i,j}\in V,\,k=0,\ldots,N-1\}
\end{equation*}
 a cylinder where the values of $\sigma_{i,j,k}\in\{-1,+1\}$ are assigned on
  $V\subset\la$.  If $ V_1 $ and $V_2 $ are two disjoint sets of
  vertices, we sometime denote by $ ( \sigma_{V_1} , \sigma_{V_2} ) $
  the cylinder $\sigma_{V_1 \cup V_2}$.



We now give an inequality that will be useful in Theorem
\ref{2piani}. If $\sigma_{X}$ and $\bar{\sigma}_{X}$ are
finite-dimensional cylinders with $(\bar{\sigma}_{X})_u \geq (
\sigma_{X})_u$ for every vertex $u \in X$, then the following
relations hold (see \cite{ligget}):
\begin{equation}\label{confronti}
\mu_{\la}(\bi\,|\,\sigma_{X})\leq\mu_{\la}(\bi\,|\,\bar{\sigma}_{X}),
\hbox{ } \hbox{ } \hbox{ } \hbox{ }
\mu_{\la}(\bii\,|\,\sigma_{X})\geq\mu_{\la}(\bii\,|\,\bar{\sigma}_{X}).
\end{equation}

\begin{teo}\label{2piani}
 For a ferromagnetic
Ising model on $(\Sigma^{(2)},\F,\mu_{\pm})$ at zero external
field the inequality $ |M(\mu_{\pm})|\leq
R(\mathbf{c}^{\pm};\,\mu_{\pm}) $ holds. Moreover
$R(\mathbf{c}^{\pm};\,\mu_{\pm})>0$ if and only if $
|M(\mu_{\pm})|>0 $.
\end{teo}
\begin{proof}
 For a fixed
$Y\in\Gamma$, take a N-box $  \Lambda_0 \supset Y \cup \partial Y
$.
Let consider the cylinder $\sigma_{Y}\supset\cyi$ that we will
denote with ${\mathbf 1}_{Y}$ that assigne $+1$ spins to all the
2-vertices belonging to $Y$.  There exists also a family $ \{
\sigma_{\partial Y} \} $ of cylinders $\sigma_{\partial
Y}\subset\partial\cyi$, moreover the cylinders $\{ ({\mathbf 1}_Y,
\sigma_{\partial Y} )\}_{ \sigma_{\partial Y} \subset\partial\cyi}
$ on the vertices $(Y \cup
\partial Y )$ are disjoint and $ \bigcup_{ \sigma_{\partial Y} \subset
\partial\cyi } ({\mathbf 1}_Y, \sigma_{\partial Y} ) = \cyi $.

Then we can write
%
%
\begin{eqnarray*}
\mlo^{+}(\cyi)&=&\frac{\mu_{\la}(\cyi\cap\bi)}{\mu_{\la}(\bi)}=\frac{1}{\mu_{\la}(\bi)}\sum_{\sigma_{\partial
Y}\subset\partial\cyi}
\mu_{\la}(({\mathbf 1}_{Y},\sigma_{\partial Y})\cap B^{+})=\\
&=&\frac{1}{\mu_{\la}(\bi)}\sum_{\sigma_{\partial
Y}\subset\partial\cyi}\mu_{\la}(B^{+}\,|\,({\mathbf 1}_{Y},
\sigma_{\partial Y}))\mu_{\la}(({\mathbf 1}_{Y},\sigma_{\partial Y}))=\\
&=& \frac{1}{\mu_{\la}(\bi)}\sum_{\sigma_{\partial
Y}\subset\partial\cyi}\mu_{\la}(B^{+}\,|\,\sigma_{\partial Y})\,
\mu_{\la}({\mathbf 1}_{Y} \,|\,\sigma_{\partial
Y})\,\mu_{\la}(\sigma_{\partial Y}) .
\end{eqnarray*}
where we are using Markov property in the last equality.

Similarly
\begin{equation*}
\mlo^{-}(\cyi)=  \frac{1}{\mu_{\la}(\bii)}\sum_{\sigma_{\partial
Y}\subset\partial\cyi}\mu_{\la}(B^{-}\,|\,\sigma_{\partial Y})\,
\mu_{\la}({\mathbf 1}_{Y} \,|\,\sigma_{\partial
Y})\,\mu_{\la}(\sigma_{\partial Y}) .
\end{equation*}
Since $\mu_{\la}(B^{+})=\mu_{\la}(B^{-})$, we have
\begin{equation}\label{sup}
\frac{\mlo^{+}(\cyi)}{\mlo^{-}(\cyi)} =
\frac{\sum_{\sigma_{\partial
Y}\subset\partial\cyi}\mu_{\la}(B^{+}\,|\,\sigma_{\partial Y})\,
\mu_{\la}({\mathbf 1}_{Y} \,|\,\sigma_{\partial
Y})\,\mu_{\la}(\sigma_{\partial Y}) }{\sum_{\sigma_{\partial
Y}\subset\partial\cyi}\mu_{\la}(B^{-}\,|\,\sigma_{\partial Y})\,
\mu_{\la}({\mathbf 1}_{Y} \,|\,\sigma_{\partial
Y})\,\mu_{\la}(\sigma_{\partial Y}) } \leq
\sup_{\sigma_{\partial
Y}\subset\partial\cyi}\frac{\mu_{\la}(\bi\,|\,\sigma_{\partial
Y})}{\mu_{\la}(B^{-}\,|\,\sigma_{\partial Y})} .
\end{equation}
Let us define the set
\begin{equation} \label{nF}
  F_{\partial Y} := \{ {\tilde{ \sigma }} \in \Sigma^{(2)} : ({\tilde{
\sigma}}_{i,j, 0} , {\tilde{ \sigma}}_{i,j, 1} ) \in L \} \subset
\partial C^+_Y ,
\end{equation}
where $ L = \{ (-1,1 ) , (1, -1) \}$ (we are not considering
$(-1,-1)$).  Using (\ref{confronti}) it is clear that the supremum
in (\ref{sup}) is achieved for $ \sigma_{\partial Y} \subset
F_{\partial Y} $. Let us define the operator $R : \Sigma \to
\Sigma $  as:
$$
(R \sigma)_{i, j , 1} = \sigma_{i, j , 0 } \hbox{  } \hbox{ and }
\hbox{  }
 (R \sigma)_{i, j , 0} = \sigma_{i, j , 1 }.
$$
The following equality is clear
\begin{equation} \label{rimane}
\mu_{\la} ( B^+ | \sigma_{\partial Y} ) = \mu_{\la} ( B^+ |(
R\sigma)_{\partial Y} )
\end{equation}
 because the first and second slab play the same role in the Ising
measure.  Moreover if $ \sigma_{\partial Y} \subset F_{\partial Y}
$ then $ ( R\sigma)_{\partial Y} = - \sigma_{\partial Y} $, and in
general is  $ \mu_{\Lambda} (B^+ | \sigma_{\partial Y} ) =
\mu_{\Lambda} (B^- | - \sigma_{\partial Y} ) $.

Thus for each $ \sigma_{\partial Y} \subset F_{\partial Y} $
\begin{equation}  \label{ovv}
 \mu_{\la} ( B^+ | \sigma_{\partial Y} ) =  \mu_{\la} ( B^+ |(
R\sigma)_{\partial Y} )= \mu_{\la} ( B^+ | -\sigma_{\partial Y}
)=\mu_{\la} ( B^- | \sigma_{\partial Y} ) .
\end{equation}
Hence, by (\ref{ovv}) and  previous argument
\begin{equation}
\label{uno} \sup_{\sigma_{\partial
Y}\subset\partial\cyi}\frac{\mu_{\la}(\bi\,|\,\sigma_{\partial
Y})}{\mu_{\la}(B^{-}\,|\,\sigma_{\partial Y})}=
\sup_{\sigma_{\partial Y}\subset F_{\partial
Y}}\frac{\mu_{\la}(\bi\,|\,\sigma_{\partial
Y})}{\mu_{\la}(B^{-}\,|\,\sigma_{\partial Y})}=             1.
\end{equation}
\comment{
\medskip

 For
all $\sigma_{Y}\subset\partial \cyi$ there exists a cylinder
$\bar{\sigma}_{\partial Y}$ which fix spins $\sigma_{i,j,0}=+1$
and $\sigma_{i,j,1}=-1$ or $\sigma_{i,j,0}=-1$ and
$\sigma_{i,j,1}=+1$ on each 2-vertex $\mathbf{c}_{i,j}\in\partial
Y$, such that $\bar{\sigma}_{\partial Y}\geq\sigma_{\partial Y}$
on $\partial Y$. By (\ref{confronti}) follows
\begin{eqnarray*}
& &
\mu_{\la}(\bi\,|\,\sigma_{\partial Y})\leq\mu_{\la}(\bi\,|\,\bar{\sigma}_{\partial Y}),\\
& & \mu_{\la}(\bii\,|\,\sigma_{\partial
Y})\geq\mu_{\la}(\bii\,|\,\bar{\sigma}_{\partial Y}),
\end{eqnarray*}
thus the supremum in (\ref{sup}) is assumed for some cylinder as
$\bar{\sigma}_{\partial Y}$.
  Then, by the configurational
symmetry and  the invariance of $\mu_{\la}$ under inversion of the
two lattices $\zii$, we have
\begin{equation}\label{uno}
\sup_{\sigma_{\partial
Y}\subset\partial\cyi}\frac{\mu_{\la}(\bi\,|\,\sigma_{\partial
Y})}{\mu_{\la}(B^{-}\,|\,\sigma_{\partial Y})}=1.
\end{equation}
By (\ref{sup}) and (\ref{uno}) follows
\begin{equation}\label{m+-}
\mlo^{+}(\cyi)\leq\mlo^{-}(\cyi).
\end{equation}
} Since this relation holds for all $\Lambda_{o}\supset
Y\cup\partial Y$ then also in the limit
$\Lambda_{o}\to\zii\times\{0,\,1\}$, by (\ref{sup}) and
(\ref{uno}) we obtain (\ref{impl.npiani}).

 To prove the inequality between 2-percolation probability and magnetization
it is enough to observe that,  by symmetry, we have
 $\E_{\mu_{\pm}}(\sigma_{0,0,0})=\E_{\mu_{\pm}}(\sigma_{0,0,1})$,
 hence
\begin{eqnarray*}
M(\mu_{\pm})& =
&\frac{1}{2}\,\E_{\mu_{\pm}}(\sigma_{0,0,0}+\sigma_{0,0,1})=\\
&=&\sum_{\overline{\sigma}_{0},\overline{\sigma}_{1}\in\{-1,\,+1\}}\frac{1}{2}\,(\overline{\sigma}_{0}+
\overline{\sigma}_{1})\,\mu_{\pm}(\sigma_{0,0,0}=
\overline{\sigma}_{0},\,\sigma_{0,0,1}=\overline{\sigma}_{1})=\\
& = & \mu_{\pm}(E^{+})-\mu_{\pm}(E^{-}).
\end{eqnarray*}
Now, the first claim of the theorem  immediately  follows by
(\ref{e+-}) and (\ref{impl.npiani}). \comment{ We also have, as in
\cite{cnpr}\setlength\arraycolsep{2pt}
\begin{eqnarray*}
\mu_{+}(E^{+}) &\!= &\!
\mu_{+}(\{\sigma\in\Sigma^{(2)}:\mathbf{c}_{0,0}\in\textrm{
($\mathbf{c}^{+}$)-finite cluster of $\sigma$}\})+\\ &\! + &\!
\mu_{+}(\{\sigma\in\Sigma^{(2)}:\mathbf{c}_{0,0}\in\textrm{
($\infty,\mathbf{c}^{+}$)-cluster of $\sigma$}\})\!\!=\\ &\! = &\!
\sum_{Y\in\Gamma}\mu_{+}(\cyi)+\mu_{+}(C_{\infty}^{+}).
\end{eqnarray*}
By (\ref{impl.npiani}) follows \setlength\arraycolsep{2pt}
\begin{eqnarray*}
M(\mu_{+})& = &
\mu_{+}(E^{+})-\mu_{+}(E^{-})=\mu_{+}(E^{+})-\mu_{-}(E^{+})=\\
& = &
\sum_{Y\in\Gamma}(\mu_{+}(\cyi)-\mu_{-}(\cyi))+\mu_{+}(C_{\infty}^{+})-\mu_{-}(C_{\infty}^{+})\leq\\
& \leq & \mu_{+}(C_{\infty}^{+}),
\end{eqnarray*}
that is
\begin{equation*}
M(\mu_{+})\leq R(\mathbf{c}^{+};\,\mu_{+}).
\end{equation*}
Taking  $(-)$-boundary conditions, we have
\begin{equation*}
|M(\mu_{-})|\leq R(\mathbf{c}^{-};\,\mu_{-}).
\end{equation*}
To prove the  second part, we project the two slabs in an unique
lattice, $\zii$, in the following way. We assign spin +1 ($-1$) if
the corresponding 2-vertex has spins +1 ($-1$) in both its
vertices, and we assign spins $+1$ or $-1$ with probability
$\frac{1}{2}$ to remaining vertices. We note that if there exists
an infinite $(\mathbf{c^{+}})$-cluster in $\zii\times\{0,1\}$,
then there exists an infinite $(\pm)$-cluster in the new lattice.
With respect to this new lattice, we consider the measure
$\pi_{\pm}$, which is the projection on $\zii$ of $\mu_{\pm}$. If
$M(\mu_{\pm})=0$ then $M(\pi_{\pm})=0$. Similarly to Theorem
\ref{one} (ii) by using the result given in \cite{non.coesistenza}
and noting that $\pi_{\pm}$ satisfy all the required hypotheses,
follows $R({\pm};\,\pi_{\pm})=0$.
 Thus also
$R(\mathbf{c}^{\pm};\,\mu_{\pm})=0$ by the observation above.
Conversely, if $|M(\mu_{\pm})|>0$, then
$R(\mathbf{c}^{\pm};\,\mu_{\pm})>0$, by the first part.} The
second claim   follows by the first inequality and
 Proposition \ref{perc.npiani}.

\end{proof}

We present another particular case, in which we are able to prove
 (\ref{impl.npiani}), and thus to obtain
characterization of phase transition via percolation. We consider
the graph ${\tilde {\mathcal G}}_3$ having vertex set
$\zii\times\{0,1,2\}$ and edge set $\E^{2,3}\cup\E^{p}$, where
\begin{equation*}
\E^{p}=\{\langle(i,j,0),(i,j,2)\rangle:i,j\in\mathbb{Z}\}.
\end{equation*}

Consider on  ${\tilde {\mathcal G}}_3$  the ferromagnetic Ising
measures $\mu_{\pm}$, and define
\begin{equation} \label{contien}
D^+_Y = \{ \sigma \in \Sigma^{(N)} : Y \in \Gamma \hbox{ belongs
to a } ({\mathbf c}^+)-\hbox{cluster of } \sigma  \} ,
\end{equation}
so that $ C^+_Y =D^+_Y \cap \partial C^+_Y $. The event $D^+_Y $
depends only on the values of $\{ \sigma_{i,j,k}:\mathbf{c}_{i,j}
\in Y,\,k=0,\ldots,N-1\} $.

\begin{teo}\label{trepiani}
For a ferromagnetic Ising model on $(\Sigma^{(3)},\F,\mu_{\pm})$,
at zero external field $R(\mathbf{c^{\pm}};\mu_{\pm})> 0 $ if and
only if $ |M(\mu_{\pm})|> 0$.
\end{teo}
\begin{proof}
If $R(\mathbf{c^{\pm}};\mu_{\pm})> 0$, then $|M(\mu_{\pm})|> 0$ by
Proposition \ref{perc.npiani} with $N=3$. Conversely, if $ |
M(\mu_{\pm}) | > 0$, we can use Proposition \ref{prop.n+1} and
prove that (\ref{impl.npiani}) holds. Indeed, for a fixed
$Y\in\Gamma$,  consider the set of all cylinders
$\sigma_{Y}\subset D^+_Y$ and $\sigma_{\partial
Y}\subset\partial\cyi$.
Then, using Markov property, we have
\begin{eqnarray*}
& &
\mlo^{+}(\cyi)=\frac{\mu_{\la}(\cyi\cap\bi)}{\mu_{\la}(\bi)}=\\
& &=\frac{1}{\mu_{\la}(\bi)}\sum_{\sigma_{\partial
Y}\subset\partial\cyi,\sigma_{Y}\subset D^+_Y}
\mu_{\la}((\sigma_{Y},\sigma_{\partial Y})\cap B^{+})=\\
& &=\frac{1}{\mu_{\la}(\bi)} \sum_{\sigma_{\partial
Y}\subset\partial\cyi}\mu_{\la}(B^{+}\,|\,\sigma_{\partial
Y})\mu_{\la}(\sigma_{\partial Y})\sum_{\sigma_{Y}\subset
D^+_Y}\mu_{\la}(\sigma_{Y}\,|\,\sigma_{\partial
Y})    ,        
\end{eqnarray*}
and similarly for  $\mlo^{-}(\cyi)$, \comment{\begin{equation*}
\mlo^{-}(\cyi)= \frac{1}{\mu_{\la}(\bii)}\sum_{\sigma_{\partial
Y}\subset\partial\cyi}\mu_{\la}(\sigma_{\partial Y}\cap
B^{-})\sum_{\sigma_{Y}\subset\cyi}\mu_{\la}(\sigma_{Y}\,|\,\sigma_{\partial
Y}).
\end{equation*}}
hence:
\begin{eqnarray}\label{suptre}
& & \frac{\mlo^{+}(\cyi)}{\mlo^{-}(\cyi)}\leq\comment{&=
&\frac{\sum_{\sigma_{\partial
Y}\subset\partial\cyi}\mu_{\la}(\sigma_{\partial Y}\cap
B^{+})\sum_{\sigma_{Y}\subset\cyi}\mu_{\la}(\sigma_{Y}\,|\,\sigma_{\partial
Y})}{\sum_{\sigma_{\partial
Y}\subset\partial\cyi}\mu_{\la}(\sigma_{\partial Y}\cap
B^{-})\sum_{\sigma_{Y}\subset\cyi}\mu_{\la}(\sigma_{Y}\,|\,\sigma_{\partial Y})}\leq \\
 & \leq &
\sup_{\sigma_{\partial
Y}\subset\partial\cyi}\frac{\mu_{\la}(\sigma_{\partial
Y}\cap\bi)}{\mu_{\la}(\sigma_{\partial Y}\cap
B^{-})}=}\sup_{\sigma_{\partial
Y}\subset\partial\cyi}\frac{\mu_{\la}(\bi\,|\,\sigma_{\partial
Y})}{\mu_{\la}(B^{-}\,|\,\sigma_{\partial Y})}\nonumber.
\end{eqnarray}
We need to prove that
\begin{equation}\label{minoreuno}
\sup_{\sigma_{\partial
Y}\subset\partial\cyi}\frac{\mu_{\la}(\bi\,|\,\sigma_{\partial
Y})}{\mu_{\la}(B^{-}\,|\,\sigma_{\partial Y})}\leq 1.
\end{equation}
Let us define
\begin{equation}\label{rit}
G_{\partial Y } : = \{ {\tilde{ \sigma }} \in \Sigma : ({\tilde{
\sigma}}_{i,j, 0} , {\tilde{ \sigma}}_{i,j, 1} , {\tilde{
\sigma}}_{i,j, 2} ) \in L \} \subset \partial C^+_Y ,
\end{equation}
where $L = \{ (-1,-1,1), (-1,1,-1), (1,-1,-1),   \}$. Relations
(\ref{confronti}) shows that the supremum in (\ref{minoreuno}) is
achieved on cylinders that are subset of $ G_{\partial Y } $.

For total spin flip invariance
\begin{equation}\label{relunop}
\mu_{\la}(\bi\,|\, \sigma_{\partial Y}) = \mu_{\la}( B^-\,|\, -
\sigma_{\partial Y}) .
\end{equation}
We now define the rotation operator $ R : \Sigma \to \Sigma $ as:
\begin{equation}\label{rotaz}
(R \sigma)_{i,j,k} = \sigma_{i,j,k+1} \hbox{ } \hbox{ } \hbox{ }
\hbox{ } \forall (i,j,k ) \in {\mathbb Z}^2 \times \{ 0, 1,2
\},k=0,1,2
\end{equation}
where $ \sigma_{i,j,0}= \sigma_{i,j,3}$. Since $\mu_{\la}$ is
invariant under rotation of the three slabs, we have
\begin{equation}\label{relno}
\mu_{\la}(\bi\,|\,\sigma_{\partial Y})=\mu_{\la}(\bi\,|\,( R
\sigma)_{\partial Y}) = \mu_{\la}( B^- \,|\,  (- R
\sigma)_{\partial Y}).
\end{equation}
Now  observe that if $ \mathbf{c}_{i,j} \in \partial Y $,  then $
\sigma_{i,j,k } \leq (- R \sigma)_{ i,j,k } $ holds for all
$\sigma \in G_{
\partial Y}$, hence
\begin{equation}\label{runol}
\mu_{\la}(\bi\,|\,\sigma_{\partial Y})= \mu_{\la}( B^- \,|\,  (- R
\sigma)_{\partial Y}) \leq
 \mu_{\la}( B^- \,|\,  \sigma_{\partial Y}) ,
\end{equation}
and
\begin{equation}\label{glg}
 \limsup_{\Lambda \uparrow\zii\times\{0,1,2\}}   \sup_{\sigma_{\partial
Y}\subset\partial\cyi}\frac{\mu_{\la}(\bi\,|\,\sigma_{\partial
Y})}{\mu_{\la}(B^{-}\,|\,\sigma_{\partial Y})}\leq 1,
\end{equation}
implying  (\ref{impl.npiani}). This concludes the proof.
\end{proof}
Theorem \ref{trepiani} says that there exists a phase transition
in the Ising model on $\Sigma^{(3)}$ if and only if there is a
positive  probability that the 3-vertex at the origin belongs to
an infinite cluster of 3-vertices with a majority of $+1$ spins on
its vertices. Contrary to the case of two slabs, we do not obtain
an inequality between the 3-percolation probability and
magnetization since in the case of three slabs we cannot write
\begin{equation*}
M(\mu_{\pm})=\mu_{\pm}(E^{+})-\mu_{\pm}(E^{-}).
\end{equation*}

A natural problem to address is that of determining whether there
exists a maximal number of slabs for which the only extremal Gibbs
measures are $\mu_+$ and $\mu_-$.
\begin{problem}\label{finale}
Define
$$
N_{c} = \sup \{ N \in {\mathbb N}: \hbox{ the only extremal
measures on $N$-slabs are $\mu_+$ and $\mu_-$} \}  .
$$
  Two natural questions are: is 
$N_c $ finite or infinite?  Is $N_c$ equal to one?  We
conjecture that $N_c= \infty$. Hence the behavior in three
dimensions should remain meaningfully different with respect to
slab graphs.
\end{problem}


\medskip

We end the paper with some remarks.

\medskip

Throughout the paper we have only considered  constant
interactions equal to $1$.  However, one can see that all  proofs
work similarly if  different values of the interactions on the
slabs are chosen: $J_o $ interaction between spins on the same
level and $J_v $ vertical interaction between spins on different
slabs. The symmetries between the slabs in Theorem~\ref{2piani}
and Theorem~\ref{trepiani} still hold and so the proofs work
without modifications.

\medskip

A second remark concerns a possible interpretation  of Proposition
\ref{prop.n+1}. Indeed, Proposition \ref{prop.n+1} can be used in
order to obtain a filtering result; let us suppose that a
configuration $\sigma\in\Sigma^{(N)} $ generated by the measure
$\mu_+ $ is represented only by giving  the following information:
on the $N$-vertex $v $ there is a $ +1$ majority, a $-1$ majority
or the same proportion of $+1$ and $-1$. By using Proposition
\ref{prop.n+1} we can say that this information is sufficient to
establish whether $\mu_{+}$ is in a region of phase transition ($T
< T_c $) or not ($T \geq T_c $). It's enough to observe that, on a
sequence of boxes invading all the space, the frequency of
$N$-vertices with $+1$ majority is definitively larger than the
frequency of $N$-vertices with $-1$ majority if and only if there
is phase transition (we are also using the ergodicity of the
measure $\mu_+$).

\comment{

Una seconda osservazione riguarda la Propsizione \ref{prop.n+1} ed
una sua interpretazione. Possiamo pensare ad un problema di
filtraggio in cui una configurazione sugli N-piani della misura
$\mu_+$ viene rappresentata dando solo l'informazione se su un
N-vertice c'è maggioranza di $+1$ (oppure di $-1$). Dalla
Proposizione \ref{prop.n+1} si evince che questa informazione è
suffiicnte per determinare se siamo in un regime di transizione di
fase o no. Infatti basta solo osservare, su regioni sempre più
ampie, e vedere se la frequenza degli N-vertici con maggioranza di
$+1$ è diversa da quella degli N-vertici con maggioranza di $-1$.
Il sistema essendo ergodico ammette il limite (per $\Lambda \to
\mathbb Z ^2 \times \{ 0 ,1, \dots , N-1 \} $) e il limite stesso
sarà quasi certamente costante. }

\medskip

We present also an extension of the phase transition
characterization via percolation to some exotic graphs. We only
give an example of these graphs in which the result can be
applied. Let's consider  $\mathbb L ^2 =(\mathbb Z ^2, \mathbb E
^2 )$.  For each vertex $v=(i,j)  \in \mathbb Z ^2$ we take a
number $n_0 $ of vertices  denoted by $(v;l)=(i,j; l)$ for $l = 1,
\dots , n_0 $.  The set $H_{v}=\{(v;l),l = 1, \dots , n_0\}$ is
called \emph{hyper-vertex}.  We put an edge between all pairs of
vertices
$(v;l),(v;m),l,m =1, \dots , n_0$; 
moreover for $e =<u,v> \in \mathbb E ^2 $ we set an edge between
each pair of vertices $(u;l),(v;m),l,m =1, \dots , n_0$.
 Now let's
define the Ising model with plus boundary conditions on this graph
and declare that the \emph{hyper-spin} $S_v $ on the hyper-vertex
$H_v$ is equal to $ \hbox{sign} (\sum_ {l = 1 }^{ n_0 } \sigma_{
v; l }) $, where $\hbox{sign} (0) =0 $.
 For
the random field $ \{ S_v  \}_{ v \in {\mathbb Z}^2 } $ there is
percolation (\emph{i.e.} there exists an infinite cluster of
hyper-vertices with plus hyper-spins) if and only if $\mu_+$ is in
the phase transition region ($T < T_c$).

\comment{
 Now put edges and interactions between all the pairs of
vertices $ \{ (v;l):{ l=1, \dots , n_0} \}$ (they form a complete
graph), also consider the following: if $ e =<u,v> \in \mathbb E
^2 $ then are edges $ \{ <(u;l),(v;m)> :{ l,m =1, \dots , n_0}
\}$. For the graph }

\bibliography{slab6giu}

\end{document}